\documentclass[10pt]{article}
\usepackage[utf8]{inputenc}
\usepackage{latexsym,amssymb,amsmath,url,authblk}

\def\N{\mathbb{N}}

\def\Z{\mathbb{Z}}
\def\R{\mathbb{R}}

\def\proof{\par\noindent{\em Proof. }}
\def\eproof{\hfill{$\Box$}\bigskip}

\def\tec{\hspace{-1.6mm}{\bf. }}
\def\ds{\dots}
\def\sus{\subset}
\def\al{\alpha}
\def\be{\beta}
\def\ga{\gamma}

\def\cc{\colon}
\def\ep{\varepsilon}

\newtheorem{thm}{Theorem}[section]
\newtheorem{prop}[thm]{Proposition}
\newtheorem{cor}[thm]{Corollary}

\title{Bertrand's paradox on a~monitor}
\author{Martin Klazar\footnote{klazar@kam.mff.cuni.cz}}
\date{\today}

\begin{document}

\maketitle

\begin{abstract}
We investigate Bertrand's probabilistic paradox through the lens of discrete geometry 
and old-fashioned but reliable discrete probability. We approximate the plane 
unit circle with $\frac{1}{n}\times\frac{1}{n}$ boxes and count the pairs 
of boxes separated by distance more than $\sqrt{3}$. For $n\to\infty$ the proportion 
of such pairs goes to 
$$
\frac{1+\sqrt{3}}{8}-\frac{\pi(2-\sqrt{3})}{96}=0.33273\ds\;.
$$
\end{abstract}

\section{Introduction}

We discuss and implement a~discrete approach to the well 
known Bertrand's probabilistic paradox. J.~Bertrand proposed it 
in his textbook \cite[pp.~5--6]{bert} in 1889. Another, less known, Bertrand's paradox 
is in economics \cite{buch_felt,gero,hehe}.
The probabilistic paradox requires to compute the probability 
$$
\mathrm{Pr_{Be}}
$$
that a~random chord $T$ in a~plane circle $C$ is longer than a~side of the inscribed equilateral triangle. If $C$ has radius $1$, then $T$ 
has to be longer than $\sqrt{3}$. J.~Bertrand proposed three numerically 
different solutions, given below as Solutions 1--3. It is paradoxical that the solution is far from unique. 

{\em Solution~1. }For $T$ determined by its endpoints $\{a,b\}=C\cap T$ one gets 
$\mathrm{Pr_{Be}}=\frac{1}{3}$.
{\em Solution~2. }Parametrization of $T$ by the intersection point of $T$ with the diameter 
$P$ of $C$ perpendicular to $T$ leads to $\mathrm{Pr_{Be}}=\frac{1}{2}$. 
{\em Solution~3. }Parametrization of $T$ by its midpoint yields $\mathrm{Pr_{Be}}=\frac{1}{4}$. 

But one does not have to stop here. {\em Solution~4. }W.~W.~Funkenbusch \cite{funk}
fixes a~point $a\in C$, selects at random another point $b$ on or inside $C$, and takes 
the chord $T$ going through $a$ and $b$. By computing the corresponding area inside $C$, 
he gets that
$$
\mathrm{Pr_{Be}}=\frac{1}{3}+\frac{\sqrt{3}}{2\pi}=0.60899\ds\;.
$$ 
This constant appears in a~different context also in \cite{holb_kim}. 
{\em Solution~5. }Inspired by the previous solution, we propose to fix a~point 
$a\in C$, select uniformly and at random a~number $d\in(0,2]$, flip a~fair coin $c\in\{h,t\}$, 
and get $T$ as the upper or lower chord, depending on the coin flip, with the endpoint $a$ 
and length $d$. Thus
$$
\mathrm{Pr_{Be}}=2-\sqrt{3}=0.26794\ds\;.
$$
{\em Solution~6. }O.~Sheynin \cite{shey} mentions that R.~De Montessus \cite{demo} proposed for 
$\mathrm{Pr_{Be}}$ uncountably many values. 

Let us discuss the first five solutions more abstractly. We denote the open 
interior disc of $C$ by $D$, its closure by $\overline{D}=D\cup C$, the center 
of $D$ by $s$ and we set
$$
S=C\times C\setminus\{(a,\,a)\;|\;a\in C\}\;.
$$
We denote by 
$$
\mathcal{CH}=\{\{a,\,b\}\;|\;a,\,b\in C,\,a\ne b\}
$$ 
the set of all chords of $C$, and by $\mathcal{CH}_l\subset\mathcal{CH}$ the set of chords longer than $\sqrt{3}$.
A~unified approach to the above solutions selects a~probability space $(\Omega,\Sigma,\mu)$, an (almost) surjective map 
$$
F\cc\Omega\to\mathcal{CH}
$$ 
such that $F^{-1}[\mathcal{CH}_l]\in\Sigma$, and then sets
$$
\mathrm{Pr_{Be}}=\mu(F^{-1}[\mathcal{CH}_l])\;.
$$
The following spaces are employed; by the ``product probability measure'' we mean product of 
normalized Lebesgue measures. 

Solution~1 uses $S$ with the product probability measure. 
Solution~2 uses $[0,\pi)\times(0,2)$ with the product probability measure.  
Solution~3 uses $D\setminus\{s\}$ with the normalized two-dimensional Lebesgue measure. 
Solution~4 uses
$$
(C\times\overline{D})\setminus\{(a,a)\;|\;a\in C\}
$$ with the product probability measure. In Solution~5 we use the space 
$$
C\times(0,2]\times\{h,t\}
$$ with the product of normalized Lebesgue measures and the counting measure on $\{h,t\}$.

In Solutions~2 and~3 the map $F$ is injective. In Solution~1 the map $F$ is two-to-one, 
two pairs $(a,b)$ always determine one chord. In Solution~5 the map $F$ is two-to-one, 
except for $d=2$ when four triples $(a,d,c)$ determine a~single chord. In Solution~4 
the map $F$ is many-to-one. No exposition of Bertrand's paradox we know of notices that 
Solution~1 differs from Solutions~2 and~3 in the non-injectivity of $F$. True, this point 
is not very important because one may replace the onto and two-to-one map 
$F\colon S\to\mathcal{CH}$ with a~bijection $G\colon S'\to\mathcal{CH}$, where $S'$ 
is a~half of $S$; this replacement preserves the value $\mathrm{Pr_{Be}}=\frac{1}{3}$ because 
the measure is normalized to give $S'$ probability $1$.

In Solutions~1, 2, 4 and~5 the map $F$ is onto but in Solution~3 the image $F[\Omega]$ misses the chords $T$ containing the center $s$. 
One may be tempted to dismiss them on the ground that they have probability~$0$, 
but this is not rigorous. These chords are outside the image of $F$ and, in the chosen model, simply have no probability. It is formally more satisfying to extend $F$ 
in the obvious way to a~bijection
$$
F\colon\Omega:=(D\setminus\{s\})\cup u\to\mathcal{CH}\;,
$$
where $u=[0,\pi)\subset\R^2$ is taken as a~plane segment disjoint to $D$, and to use in $\Omega$ again the normalized two-dimensional 
Lebesgue measure. Solution~4 shows that the tacit assumption of (almost) injectivity of $F$ is in reality not very strongly justified. One produces 
chords $T$ in $C$ by random experiments, 
and it may well happen that many of them yield the same outcome. 

We give one more {\em Solution~7. }Since both sets $\mathcal{CH}_l$ 
and $\mathcal{CH}\setminus\mathcal{CH}_l$ have the same cardinality $\mathfrak{c}$ of continuum, there is 
a~bijection $F\colon[0,1]\to\mathcal{CH}$ such that 
$F^{-1}[\mathcal{CH}_l]=[0,\frac{2021}{2022}]$. So, using the normalized 
Lebesgue measure on $[0,1]$, we get
$$
\mathrm{Pr_{Be}}=\frac{2021}{2022}=0.99950\ds\;.
$$

There is an extensive debate on Bertrand's paradox: 
\cite{aert_sass}, \cite[Chapter~1.3]{ande}, \cite{bang}, \cite{cast_dolc}, \cite{funk}, \cite{gyen_rede}, \cite{holb_kim}, \cite{kaus}, \cite{mari}, 
\cite{demo}, \cite[pp. 62--64]{renyi} (here this author learned about the paradox for the first time), \cite{rizz}, 
\cite{rowb}, \cite{rowb_shac}, \cite{shac}, \cite{shey}, \cite{sora_volc}, \cite{vido}, and \cite{wang_jack}. Most 
of this list was compiled by searching in the {\em Mathematical Reviews} for articles having 
``Bertrand'' and ``paradox'' in the title. Many more references could be added to it.

We want to contribute to this debate by proposing a~realistic mathematical model of the paradox when it is treated as 
a~physical experiment. 
The idea of experimental resolution of the paradox is an old one, but in our opinion physical 
schemes proposed so far do not  employ faithful mathematical models. For example, \cite{holb_kim} 
suggests to 
perform the random selection of $T$ by letting cosmic rays pass through a~disc chamber detector, or by imagining 
them doing so. But the rays are in \cite{holb_kim} treated as lines with zero width. 
The problem with such approaches is that nothing like a~ray or a~circle with zero width exists in reality. 
Physical chords and physical circles have nonzero, positive widths and are made of discrete parts with positive sizes. 
Any ``simplification'' that replaces a~positive number with zero 
is a~jump in a~(totally) different world. In each case it should be  carefully justified and one should check that contact with 
the original setup was not lost and that after the simplification one is not solving a~(completely) different problem.
\begin{quote}
We think that the plethora of solutions for Bertrand's paradox demonstrate that by working with chords and circles with zero width 
we simplify too much.
\end{quote}
We return to what it means, or should mean\,---\,in our opinion\,---\,to resolve Bertrand's paradox in the last section.

So we want to compute the probability in Bertrand's paradox in the old-fashioned but safe way 
as the limit ratio of two {\em finite} numbers, namely that of favorable cases and that of all cases. We denote 
this discrete version of $\mathrm{Pr_{Be}}$ by 
$$
\mathrm{Pr_{Be,\,\Box}}\;.
$$ 
We work in the Euclidean plane $\R^2$ and  for $n\to\infty$ consider a~sequence of square grids of axes-parallel 
$\frac{1}{n}\times\frac{1}{n}$ squares, so called boxes, with a~vertex of the grid in the origin. We consider physical 
unit circles $C(n)$ consisting of boxes whose interiors intersect the geometric, origin-centered unit circle 
$C$. Let $\mathrm{Ch}(n)\subset C(n)\times C(n)$ be the pairs $(B,B')$ of boxes 
in $C(n)$ with distance $>\sqrt{3}$, i.e., such that some points $b$ and $b'$ in the interior of $B$ 
and $B'$, respectively, have distance $|bb'|>\sqrt{3}$. We denote the cardinality of a~finite set $X$
by $|X|$.

In Theorem~\ref{thm_BrSquare} we compute that
\begin{eqnarray*}
  \mathrm{Pr_{Be,\,\Box}}&=&\lim_{n\to\infty}\frac{|\mathrm{Ch}(n)|}{|C(n)\times C(n)|-|C(n)|}=\lim_{n\to\infty}\frac{|\mathrm{Ch}(n)|}{|C(n)|^2}\\
  &=&\frac{1+\sqrt{3}}{8}-\frac{\pi(2-\sqrt{3})}{96}=0.33273\ds\;.
\end{eqnarray*}
The chords in $\mathcal{CH}$ are unordered pairs while in $C(n)\times C(n)$ they are ordered pairs, but the factor $2$ cancels in the ratio.
In Corollary~\ref{cor_Andel}  we determine the limit proportion
$$
\mathrm{Pr_{1,\,3,\,\Box}}:=\lim_{n\to\infty}\frac{|C_{0,\,\pi/3}(n)|}{|C(n)|}=\frac{1+\sqrt{3}}{16}=0.17075\ds
$$
of the boxes lying in the first quadrant of $C(n)$ in the angular positions $(0,\pi/3)$. In Corollary~\ref{cor_sol3} we note that there are many more (ordered or unordered) pairs of boxes in $C(n)$ than there are boxes intersecting the disc $D$. 
Thus for square grids Bertrand's Solution~3 cannot be implemented, there are not enough midpoints to represent bijectively all chords. Corollary~\ref{cor_expl_form} gives the  formula
$$
|C(n)|=8n-4\sum_{\substack{d\,|\,n^2\\ d=2e+1}}(-1)^e
$$
for the size of the physical unit circle.
For example, $|C(3)|=24-4(1-1+1)=20$. 
But\,---\,despite its appearance\,---\,the formula is not efficient; we do not know how to factorize integers efficiently, and therefore the formula does 
not tell us how to compute for given $n$ in time polynomial in $\log n$ the value $|C(n)|$. 
Propositions~\ref{prop_onExceptionals}--\ref{prop_aproximacni1} are auxiliary. In Section~4 we give some concluding remarks and address suggestions of one referee.

Circles seen today are often sets of pixels on the monitor of a~computer or  a~TV set, selected from a~fine rectangular array. Or maybe from an array with other shape, but this is not the point. Such circles inspired us to write this article, not the ideal 
geometric circle with zero width. Approximating $C$ by $C(n)$ seems inadequate for circles drawn in media like plain paper without 
the raster of graph paper, or on blackboard, slate or in sand. Is there a~realistic material such that when one draws $C$ on it, 
the distribution of mass in the approximation 
is roughly uniform with respect to the angular measure? In Bertrand's Solution~1 this is taken for 
granted. We hope to address this and similar questions at another occasion.

\section{Counting vertical and horizontal boxes}

By $\N=\{1,2,\ds\}$ we denote the set of natural numbers, $\Z$ are the integers and $\N_0=\N\cup\{0\}=\{0,1,\ds\}$. 
For $n\in\N$ we set $[n]=\{1,2,\ds,n\}$, $[0]=\emptyset$, 
for $a,b\in\Z$ with $a\le b$ we set $[a,b]_{\Z}=\{a,a+1,\ds,b\}$;   $[a,b]_{\Z}=\emptyset$ for $a>b$. For $n\in\N$ 
and $i,j\in\Z$ we consider the $\frac{1}{n}\times\frac{1}{n}$ {\em boxes}
$$
B=B(n)=B_{i,\,j}(n)=[i/n,\,(i+1)/n]\times[j/n,\,(j+1)/n]\subset\R^2\;.
$$
$B^o$ denotes the interior of the box $B$. Let $$
C=\{(x,\,y)\in\R^2\;|\;x^2+y^2=1\}
$$ 
be the geometric, origin-centered unit circle 
in the plane $\R^2$ and let $D$, $\overline{D}$ and $s=(0,0)\in D$ be as above in Section~1. We define the approximations $C(n)$ of $C$ by boxes as
$$
C(n)=\{B_{i,\,j}(n)\;|\;i,\,j\,\in\Z\wedge B_{i,\,j}(n)^o\cap C\ne\emptyset\},\ n\in\N\;.
$$
More generally, for any $n\in\N$ and any $X\sus\R^2$ we define
$$
X(n)=\{B(n)\;|\;B(n)^o\cap X\ne\emptyset\}\;.
$$
We say that a~box is {\em exceptional} if one of its four corners lies on $C$.

\begin{prop}\tec\label{prop_onExceptionals}
For any fixed $\varepsilon>0$ and all $n\in\N$, there are at most $O_{\varepsilon}(n^{\varepsilon})$ exceptional boxes.
\end{prop}
\proof
$B_{i,j}(n)$ is exceptional if and only if one of the four sums $i^2+j^2$, $(i+1)^2+j^2$, $i^2+(j+1)^2$ 
and $(i+1)^2+(j+1)^2$ equals $n^2$. Hence the number of exceptional 
boxes is $O(r_2(n^2))$, where for $m\in\N$ the 
quantity $r_2(m)$ is the number of solutions $(i,j)\in\Z^2$ 
of the equation $i^2+j^2=m$. There is the formula 
$$
r_2(m)=4\sum_{d\,|\,m\,\wedge\, d=2e+1}(-1)^e\;,
$$
see \cite{hard_wrig}. Thus the number of exceptional boxes is $O(\tau(n^2))$, 
where $\tau(m)$ is the number of divisors of $m$. It is well known (\cite{hard_wrig}) 
that $\tau(m)=O_{\varepsilon}(m^{\varepsilon})$.
\eproof

A~{\em chord $T(n)$ in $C(n)$} is any ordered pair of distinct boxes,
$$
T(n)=(B,B')\in C(n)\times C(n)\;\text{ and }\;B\ne B'\;.
$$ 
The {\em full chord $T_{B,B'}(n)$ in $C(n)$} generated by the chord $(B,B')$ is the set
$$
T_{B,\,B'}(n)=\{A\;|\;\exists\,ab:\ a\in B^o\wedge b\in (B')^o\wedge ab\cap A^o\ne\emptyset\}\;.
$$
of all boxes $A$ with interior intersected by a~straight segment $ab\sus\R^2$ that joins two points
inside $B$ and $B'$. But we will not use full chords. The main reason is that one may have
$T_{A,A'}(n)=T_{B,B'}(n)$ for two chords $(A,A')$ and $(B,B')$ even if $\{A,A'\}\ne\{B,B'\}$. 
This happens for example when 
$$
A=B_{-n,\,0}(n),\ A'=B_{n,\,-1}(n),\ B=B_{-n,\,-1}(n)\;\text{ and }\;B'=B_{n,\,0}(n)
$$
---\,then $T_{A,A'}(n)=T_{B,B'}(n)$ is formed by the two layers of boxes intersecting the 
segment $(-1,1)\times\{0\}$. We believe that this happens infrequently, so that 
$$
\frac{\text{the number of full chords in $C(n)$}}{\text{the number of chords in $C(n)$}}\to\frac{1}{2}
$$ 
as $n\to\infty$. 

A~key tool in our computations is Abel's summation formula (\cite[Thm.~0.3]{tene}). It has a~nice and
short proof, which we include here. Our proof differs from that in \cite{tene} and we actually 
give a~more general theorem. Recall that every bounded function $f\cc[a,b]\to\R$ with only 
finitely many discontinuities is Riemann-integrable; we only use Riemann integration (\cite{tene} uses Stieltjes integrals). For a~sequence $(a_n)=(a_1,a_2,\ds)\sus\R$ 
and a~real number $x$ we set ($n\in\N$)
$$
A(x)=\sum_{n\le x}a_n\;,
$$
where the empty sum is as usual defined to be $0$. For a~function 
$f\colon[a,b]\to\R$ we use the notation $[f]_a^b=f(b)-f(a)$.

\begin{prop}[Abel's summation]\tec\label{prop_sumAbel}
Let $l,m\in\N_0$ with $l<m$, $(a_n)\sus\R$ be a~sequence and $f\colon[l,m]\to\R$ 
be a~continuously differentiable 
function. Then
$$
\sum_{l<n\le m}a_nf(n)=[A(x)f(x)]_l^m-\int_l^m A(t)f'(t)\;\mathrm{dt}\;.
$$
\end{prop}
\proof
Both expressions equated in the identity are additive. This means that if $k,l,m\in\N_0$ 
with $k<l<m$, then the value of the expression on the interval $[k,m]$ equals to the sum 
of its values on the intervals $[k,l]$ and $[l,m]$. Thus it suffices to prove the identity 
for every $l\in\N_0$ and $m=l+1$. Then the left-hand side is $a_mf(m)$ and the right-hand side is
\begin{eqnarray*}
[A(x)f(x)]_l^m-A(l)\int_l^m f'&=&[A(x)f(x)]_l^m-A(l)[f(x)]_l^m
=[A(x)]_l^mf(m)\\
&=&a_mf(m)
\end{eqnarray*}
too.
\eproof

\noindent
See \cite[p.~13/14]{titc} for another application of this technique proving sum-integral identities.

We orient $C$ counter-clockwisely and call a~box $B\in C(n)$ {\em vertical} (resp. {\em horizontal}) 
if $C$ enters $B$ in the vertical (resp. horizontal) side of $B$. Each non-exceptional box in $C(n)$ is either 
vertical or horizontal, and if a~box in $C(n)$ is both then it is exceptional. For any two angles 
$\al\le\be$ in $\R$ with $\be-\al\le 2\pi$ we define 
$$
C_{\al,\,\be}\sus C
$$ to be the open arc of points $p\in C$ such that the length of the arc $P_0p\sus C$, where $P_0:=(1,0)$, 
that goes counter-clockwisely from $P_0$ to $p$ lies, modulo $2\pi$, in the interval $(\al,\be)$. We 
consider the approximation $C_{\al,\be}(n)$ of the arc $C_{\al,\be}$ by boxes
and denote by $C_{\al,\be,\mathrm{v}}(n)$ (resp. $C_{\al,\be,\mathrm{h}}(n)$) the subset 
of vertical (resp. horizontal) boxes.

\begin{prop}\tec
For any fixed $\al,\be\in\R$ with $0\le\al\le\be\le\pi$ and all $n\in\N$, we have the estimate
$$
|C_{\al,\,\be,\,\mathrm{v}}(n)|=[-\cos t]_{\al}^{\be}\cdot n+O(1)\;.
$$
For $\pi\le\al\le\be\le 2\pi$ this estimate holds with the function $\cos t$.
For any fixed $\al,\be\in\R$ with $-\pi/2\le\al\le\be\le\pi/2$ and all $n\in\N$, we have the estimate
$$
|C_{\al,\,\be,\,\mathrm{h}}(n)|=[\sin t]_{\al}^{\be}\cdot n+O(1)\;.
$$
For $\pi/2\le\al\le\be\le 3\pi/2$ this estimate holds with the function $-\sin t$.
\end{prop}
\proof
Let an $n\in\N$ and $\al,\be\in[0,\pi]$ with $\al\le\be$ be given. Let $I\sus[-n+1,n]_{\Z}$ be the interval of the integers 
$i$ such that the line $x=i/n$ intersects the arc $C_{\al,\be}$. Since the boxes $B$ in $C_{\al,\be,\mathrm{v}}(n)$ 
1-1 correspond, except possibly for the rightmost one, to $I$ (the right side of the box $B$, where $C_{\al,\be}$ enters $B$, 
lies on the line $x=i/n$, $i\in I$), we have the bounds
$$
|I|\le|C_{\al,\,\be,\,\mathrm{v}}(n)|\le|I|+1\;.
$$
For $i\in [-n,n]_{\Z}$, 
let $\al_i\in[0,\pi]$ be the angle at the vertex $s=(0,0)$ in the triangle $p_i\,s\,(1,0)$,  where $p_i\in C$ is the intersection 
of the upper closed semicircle of $C$ with the line $x=i/n$. Clearly, $i/n=\cos\al_i$. Thus
$$
|I|=\sum_{i\in I}1=n\sum_{i\in I}(\cos\al_{i+1}-\cos\al_i)=n(\cos\al_{k+1}-\cos\al_j)\;,
$$
where $j=\min I$ and $k=\max I$. Since $\cos\al\in[k/n,(k+1)/n]$, the numbers $\cos\al$ and $\cos\al_{k+1}$ differ 
by at most $1/n$. Similarly for $\cos\be$ and $\cos\al_j$. Thus by the above displayed bounds,
$$
|C_{\al,\,\be,\,\mathrm{v}}(n)|=n(\cos\al-\cos\be)+O(1)\;,
$$
which is the first stated estimate. For the other range of $\al$ and $\be$ and for the second stated estimate for horizontal boxes we argue similarly.
\eproof

\begin{prop}\tec\label{prop_aproximacni1}
For any fixed $\ep>0$, $\al,\be\in\R$ with $0\le\al\le\be\le\pi/2$ and all $n\in\N$, we have the estimate
$$
|C_{\al,\,\be}(n)|=[\sin t-\cos t]_{\al}^{\be}\cdot n+E_n-\Delta_n\;,
$$
where 
$$
E_n=O(1)\;\text{ and }\;0\le\Delta_n\le\min\big([\sin t-\cos t]_{\al}^{\be}\cdot n+E_n,\,O_{\ep}(n^{\ep})\big)\;. 
$$
In the other three quadrants of $C$ this estimate holds with the appropriate signs of $\sin t$ and $\cos t$, 
$+$ if the function increases and $-$ if it decreases.
\end{prop}
\proof
This follows from the previous proposition and Proposition~\ref{prop_onExceptionals}: 
$[\sin t-\cos t]_{\al}^{\be}\cdot n+E_n$ is the sum of estimates of the numbers of 
vertical and horizontal boxes in $C_{\al,\,\be}(n)$ and $\Delta_n$ is the number 
of the boxes which are both vertical and horizontal (and hence are exceptional).
\eproof

\begin{cor}\tec\label{cor_expl_form}
For every $n\in\N$ the cardinality of the approximate unit circle $C(n)$ equals
$$
|C(n)|=8n-4\sum_{\substack{d\,|\,n^2\\ d=2e+1}}(-1)^e\;.
$$
\end{cor}
\proof
We count the boxes in $C(n)$ more carefully:
\begin{eqnarray*}
|C(n)|&=&|C_{\mathrm{v}}(n)\cup C_{\mathrm{h}}(n)|=|C_{\mathrm{v}}(n)|+|C_{\mathrm{h}}(n)|-|C_{\mathrm{v}}(n)\cap C_{\mathrm{h}}(n)|\\
&=&|C_{\mathrm{v}}(n)|+|C_{\mathrm{h}}(n)|-|C_{\mathrm{e}}(n)|\;,
\end{eqnarray*}
where $C_{\mathrm{v}}(n)$ (resp. $C_{\mathrm{h}}(n)$) is the set of vertical (resp. horizontal) 
boxes in $C(n)$ and $C_{\mathrm{e}}(n)$ is the set of those 
boxes in $C(n)$ in which $C$ enters (counter-clockwisely) in a~vertex. Now $|C_{\mathrm{v}}(n)|=4n$ 
because exactly two
boxes in $C_{\mathrm{v}}(n)$ correspond to each of the $2n-1$ vertical lines $x=i/n$ with $i=-n+1,-n+2,\ds,n-1$ 
(the point of entrance of $C$ 
in the box lies on the line), only one box in $C_{\mathrm{v}}(n)$ corresponds to each of the two vertical lines $x=-1$ and $x=1$, and this 
exhausts all boxes in $C_{\mathrm{v}}(n)$. By symmetry, also $|C_{\mathrm{h}}(n)|=4n$.
Clearly, $C_{\mathrm{v}}(n)\cap C_{\mathrm{h}}(n)=C_{\mathrm{e}}(n)$ and $C_{\mathrm{e}}(n)$ bijectively 
corresponds to the solutions $(i,j)\in\Z^2$ of the equation $i^2+j^2=n^2$. Thus $|C_{\mathrm{e}}(n)|=r_2(n^2)$. Using
the formula in \cite{hard_wrig} for the number $r_2(m)$ of these solutions, we obtain the stated formula.
\eproof

\noindent
As we already mentioned, at the present state of knowledge it is not known how to use this formula to compute $|C(n)|$ in time polynomial in $\log n$. 

Let $B^*$ be a~fixed box, say $B^*=B_{0,0}(n)$, which will serve as an origin.
An even more ``physical'' approximations of $C$ by boxes would be
$$
C[n]=\{B(n)\;|\;\exists\,ab:\ a\in (B^*)^o\wedge b\in B(n)^o\wedge |ab|=1\},\ n\in\N\;.
$$
We hope to investigate them in \cite{klaz_buff}. Here we confine to the approximations $C(n)$ which are easier to handle.

\begin{cor}\tec\label{cor_sol3}
For any fixed $\ep>0$ and all $n\in\N$,
$$
|C(n)\times C(n)|=64 n^2+O_{\ep}(n^{1+\ep})\;.
$$
For $D=\{(x,y)\in\R^2\;|\;x^2+y^2<1\}$ one has that
$$
|D(n)|=\pi n^2+\Delta_n\;,
$$
where $0\le \Delta_n\le |C(n)|<8n$. Thus there are many more chords (ordered or unordered) in $C(n)$ than boxes in $D(n)$.
\end{cor}
\proof
The first formula follows from Corollary~\ref{cor_expl_form} and from the bound $\tau(m)=O_{\ep}(m^{\ep})$ on the number of divisors function 
(\cite{hard_wrig}), already mentioned in the proof of Proposition~\ref{prop_onExceptionals}. To prove the second formula, 
we consider besides $D(n)$ also the set $U(n)$ of boxes $B(n)$ such that $B(n)^o\sus D$, thus $U(n)\sus D(n)$. Since 
$\bigcup U(n)\sus\overline{D}\sus\bigcup D(n)$, the area of the disc $\overline{D}$ equals $\pi$ and $D(n)\setminus U(n)=C(n)$,
we have the relations
$$
\frac{|U(n)|}{n^2}\le \pi\le\frac{|D(n)|}{n^2}\;\text{ and }\;|D(n)|-|U(n)|=|C(n)|\;.
$$
By this and Corollary~\ref{cor_expl_form} we get the second formula.
\eproof

\noindent
Bertrand's Solution~3 therefore cannot be implemented for $C(n)$ and $D(n)$, there are not enough midpoints for the chords.

In the next corollary the indices $1$ and $3$ refer to the arc between the hours $1$ and $3$ on a~standard 
round clock face.

\begin{cor}\tec\label{cor_Andel}
$$
\mathrm{Pr}_{1,\,3,\,\Box}:=\lim_{n\to\infty}\frac{|C_{0,\,\pi/3}(n)|}{|C(n)|}=\frac{1+\sqrt{3}}{16}=0.17075\ds\;.
$$
\end{cor}
\proof
We set $\ep:=1/2$ and get by Proposition~\ref{prop_aproximacni1} and by the 
formula in Corollary~\ref{cor_expl_form} that
$$
\mathrm{Pr}_{1,\,3,\,\Box}=\lim_{n\to\infty}\frac{[\sin t-\cos t]_0^{\pi/3}\cdot n+
O(\sqrt{n})}{8n+O(\sqrt{n})}=\frac{1+\sqrt{3}}{16}\;.
$$
\eproof

\noindent
This is inspired by \cite[Chapter~1.3]{ande} where exposition of Bertrand's paradox 
starts with the remark that the second hand on clocks appears between hours 
$1$ and $3$ with probability $\frac{1}{6}=0.16666\ds\;$. In the world of $C(n)$ 
clocks spend slightly more time between these two hours.

\section{Counting long physical chords in $C(n)$}

In this section we compute $\mathrm{Pr_{Be,\Box}}$ by counting the {\em long chords}
$$
\mathrm{Ch}(n)=\{(B,\,B')\in C(n)^2\;|\;\exists\,ab:\ a\in B^o\wedge b\in (B')^o\wedge |ab|>\sqrt{3}\}
$$
and then computing the limit
$$
\mathrm{Pr_{Be,\,\Box}}=\lim_{n\to\infty}\frac{|\mathrm{Ch}(n)|}{|C(n)|^2}\;.
$$
We identify $C$ with the interval $[0,2\pi)$ and for $a\in C$ define the {\em arc antipodal to $a$} as
$$
\{a\}^{\mathrm{ap}}=(a+2\pi/3,\,a+4\pi/3)\ (\mathrm{mod}\ 2\pi)\ (\sus[0,\,2\pi))\;.
$$
This definition is motivated by the equivalence that
a~chord $ab$ with $a,b\in C$ has length $|ab|>\sqrt{3}$ if 
and only if $b\in\{a\}^{\mathrm{ap}}$. Bertrand's Solution~1 that $\mathrm{Pr_{Be}}=\frac{1}{3}$ follows at once. For any $X\sus C$
we define
$$
X^{\cap}=\bigcap_{a\in X}\{a\}^{\mathrm{ap}}\;.
$$
Thus 
$$
a\in X\sus C\wedge b\in X^{\cap}\Rightarrow |ab|>\sqrt{3}\;. 
$$
For $a\in C$ and $d\in(0,\sqrt{3}]$ we define more generally
$$
\{a\}^{\mathrm{ap},\,d}=\{b\in C\;|\;|ab|>d\}
$$
to be the arc of the other endpoints $b$ of chords $ab$ with length $|ab|>d$. For any $X\sus C$ we also define
$$
X^{\cup}=\bigcup_{a\in X}\{a\}^{\mathrm{ap}}\;\text{ and }\;X^{\cup,\,d}=\bigcup_{a\in X}\{a\}^{\mathrm{ap},\,d}.
$$
Clearly, $X^{\cap}\sus X^{\cup}\sus X^{\cup,d}$ and  
$$
a\in X\sus C\wedge b\in C\wedge|ab|>d\Rightarrow b\in X^{\cup,\,d}\;. 
$$

\begin{thm}\tec\label{thm_BrSquare}
The limit proportion of long chords among all chords in the approximate unit circle $C(n)$ is
$$
\mathrm{Pr_{Be,\,\Box}}=\lim_{n\to\infty}\frac{|\mathrm{Ch}(n)|}{|C(n)|^2}=\frac{1+\sqrt{3}}{8}-\frac{\pi(2-\sqrt{3})}{96}=0.33273\ds\;.
$$
\end{thm}
\proof We split the computational proof in seven steps.

\medskip\noindent
{\bf 1.~Deducing and proving equation~(\ref{intervals}) below. }Let $n\in\N$. For $j\in[4]$ we define 
$$
\mathrm{Ch}_{j\mathrm{q}}(n)=\{(B,\,B')\in\mathrm{Ch}(n)\;|\;B\in C_{(j-1)\pi/2,\,j\pi/2}(n)\}
$$ 
to be the set of long chords $(B,B')$ with $B$ in the $j$-th quadrant of $C(n)$. The sets $\mathrm{Ch}_{j\mathrm{q}}(n)$,  
$j\in[4]$, partition $\mathrm{Ch}(n)$ and,  by symmetry, have equal cardinalities. Thus
\begin{equation}\label{DnaD1qn}
|\mathrm{Ch}(n)|=4\cdot|\mathrm{Ch}_{1\mathrm{q}}(n)|
\end{equation}
and it suffices to count $\mathrm{Ch}_{1\mathrm{q}}(n)$. We employ a~function $k=k(n)\in\N$, $n\in\N$ and $k\le n$, 
such that $k\to\infty$ with $n$ but more slowly than $n$; we specify it at the end. We assume that $n>6$ and $k<n/6$ 
and refer to it as the {\em basic assumption on $n$ and $k$}. We
divide $C$ in $12k$ equally long subarcs 
$$
C_i=C_{(i-1)\ga,\,i\ga},\;\text{ where }\;i\in[12k]\;\text{ and }\;\ga=\ga(k)=\pi/6k\;.
$$
We describe the sets $C_i^{\cap}(n)$, $i\in[3k]$, and their partitions in quadrants\,---\,these are needed because 
the estimates of Proposition~\ref{prop_aproximacni1} apply only within a~single quadrant. 

For $i\in[k]$ we set 
$$
I_1=I_1(i,\,k)=[4k+i+1,\,6k]_{\Z}\;\text{ and }\;J_1=J_1(i,\,k)=[6k+1,\,8k+i-1]_{\Z}\;.
$$
For $i\in[k+1,2k]_{\Z}$ we similarly set
$$
I_2=[4k+i+1,\,6k]_{\Z},\ J_2=[6k+1,\,9k]_{\Z}\;\text{ and }\;K_2=[9k+1,\,8k+i-1]_{\Z}\;,
$$
and for $i\in[2k+1,3k]_{\Z}$ we set
$$
I_3=[4k+i+1,\,9k]_{\Z}\;\text{ and }\;J_3=[9k+1,\,8k+i-1]_{\Z}\;.
$$
We set $K_1=K_3=\emptyset$. For $j\in[3]$ and $i\in[(j-1)k+1,jk]_{\Z}$ we set 
 $$
 L(i,k)=I_j\cup J_j\cup K_j=[4k+i+1,\,8k+i-1]_{\Z}
 $$ 
 and 
 $$
 L'(i,k)=[4k+i,\,8k+i]_{\Z}\;.
 $$
For any set $I\sus[12k]$ we denote
$$
C_I=\bigcup_{i\in I}C_i\sus C\;.
$$
Thus $C_I(n)=(\bigcup_{i\in I}C_i)(n)=\bigcup_{i\in I}C_i(n)\sus C(n)$.  It follows from the above definitions that for every $i\in[3k]$,  
\begin{equation}\label{pruAsjedn}
C_i^{\cap}(n)=C_{L(i,k)}(n)\;\text{ and }\; C_i^{\cup}(n)=C_{L'(i,k)}(n)\;.
\end{equation}
(To see this, note that if $I\sus C$ is an open interval and $\overline{I}$ is its closure, 
then $I(n)=\overline{I}(n)$.)

We claim that for all $n,k\in\N$ satisfying the basic assumption,
\begin{eqnarray}
|\mathrm{Ch}_{1\mathrm{q}}(n)|&=&\sum_{i=1}^{3k}|C_i(n)|\cdot|C_{L(i,\,k)}(n)|+O(kn)+O(n^2/k)\;.\label{intervals}
\end{eqnarray}
To prove it, we first consider chords $(B,B')$ counted in the sum of products 
on the right-hand side of equation (\ref{intervals}): 
there is an $i\in[3k]$ such that $B\in C_i(n)$ and $B'\in C_{L(i,k)}(n)$. 
Since $ C_{L(i,k)}(n)=C_i^{\cap}(n)$ (by equations (\ref{pruAsjedn})), it follows that 
$(B,B')\in\mathrm{Ch}_{1\mathrm{q}}(n)$ and these chords are counted also in the left-hand side. 
But there is an overcount, a~chord $(B,B')$ is counted in the right-hand side more 
than once iff $B\in C_i(n)\cap C_{i'}(n)$ for some $i<i'$ in $[3k]$. 
This may happen only if $i'=i+1$ because $\ga>\pi/n$ (the maximum length of 
a~circular arc contained in a~$1\times 1$ square is $\pi$). Since 
$$
|C_i(n)\cap C_{i+1}(n)|=|\{i\ga\}(n)|\le 1\;\text{ and }\;|C_{L(i,\,k)}(n)|=O(k\cdot n/k)=O(n)
$$
by Proposition~\ref{prop_aproximacni1}, 
this overcount is bounded by the term $O(kn)$.

We also have to consider chords $(B,B')\in \mathrm{Ch}_{1\mathrm{q}}(n)$ 
counted in the left-hand side of equation (\ref{intervals}) but not in the right-hand side. This means that, setting $d=\sqrt{3}-2\sqrt{2}/n$ 
(each box has diameter $\sqrt{2}/n$), there is 
an $i\in[3k]$ such that $B\in C_i(n)$ but 
\begin{eqnarray*}
B'\in C^{\cup,\,d}_i(n)\setminus C^{\cap}_i(n)&=&( C^{\cup,\,d}_i(n)
\setminus C^{\cup}_i(n))\cup(C^{\cup}_i(n)\setminus C^{\cap}_i(n))\\
&=:&X_{i,\,n}\cup Y_{i,\,n}\;.
\end{eqnarray*}
For every $i\in[3k]$ and every $n>6$,
$$
|X_{i,\,n}|=O(1)\;\text{ and }\;|Y_{i,\,n}|=O(n/k)\;.
$$
The first bound follows from Proposition~\ref{prop_aproximacni1} because 
$|X_{i,n}|$ is at most the number of boxes whose interiors intersect 
$C_i^{\cup,d}\setminus  C_i^{\cup}$, which consists of two arcs of $C$ 
with lengths $O(\sqrt{3}-d)=O(1/n)$. The second bound also follows from 
Proposition~\ref{prop_aproximacni1} because 
$|Y_{i,\,n}|$ is at most the number of boxes whose interiors intersect 
two arcs of $C$ with lengths $\ga=O(1/k)$ (see equalities \ref{pruAsjedn}). 
The first bound is subsumed in the second because $k\le n$. Again by Proposition~\ref{prop_aproximacni1}, $|C_i(n)|=O(n/k)$.  Thus the surplus 
of the left-hand side in (\ref{intervals}) over the sum of products in the right-hand 
side is compensated by the term
$$
O(k\cdot(n/k)^2)=O(n^2/k)\;.
$$
Equality (\ref{intervals}) is therefore proven.

\medskip\noindent
{\bf 2.~Symmetries $S_1(n)=S_3(n)$ and $S_2'(n)=S_2'''(n)$ 
in (\ref{symmetry}). }For $j\in[3]$ we define
\begin{eqnarray}
S_j(n)&=&\sum_{i=(j-1)k+1}^{jk}|C_i(n)|\cdot|C_{L(i,k)}(n)|\label{cardinal1}\\
&=&\sum_{i=(j-1)k+1}^{jk}|C_i(n)|\cdot\big(|C_{I_j(i,\,k)}(n)|+|C_{J_j(i,\,k)}(n)|
+|C_{K_j(i,\,k)}(n)|\big)+\notag\\
&&+\;O(n)\;,\label{cardinal}
\end{eqnarray}
so 
\begin{equation}\label{Ch1qaSj}
|\mathrm{Ch}_{1\mathrm{q}}(n)|=\sum_{j=1}^{3}S_j(n)+O(kn)+O(n^2/k)\;.
\end{equation}
The error term $O(n)=O(k\cdot n/k)$ on the right-hand side of equation (\ref{cardinal})\,---\,recall that $|C_i(n)|=O(n/k)$\,---\,compensates like earlier its overcount of boxes in $(\ds+\ds+\ds)$ which occurs even though $I_j$, $J_j$ and $K_j$ are disjoint. We define
\begin{eqnarray}
S_2'(n)&=& \sum_{i=k+1}^{2k}|C_i(n)|\cdot|C_{I_2(i,\,k)}(n)|\label{S2'}\;,\\
S_2''(n)&=& \sum_{i=k+1}^{2k}|C_i(n)|\cdot|C_{J_2(i,\,k)}(n)|\;\text{ and }\label{S2''}\\
S_2'''(n)&=& \sum_{i=k+1}^{2k}|C_i(n)|\cdot|C_{K_2(i,\,k)}(n)|\label{S2'''}\;.
\end{eqnarray}
Hence 
\begin{equation}\label{S2_sumofthree}
    S_2(n)=S_2'(n)+S_2''(n)+S_2'''(n)+O(n)\;.
\end{equation}

We claim that for every $n\in\N$,
\begin{equation}\label{symmetry}
    S_1(n)=S_3(n)\;\text{ and }S_2'(n)=S_2'''(n)\;.
\end{equation}
This follows from the definition of these sums in equations (\ref{cardinal1}), 
(\ref{S2'}) and (\ref{S2'''}) and from the involution $(x,y)\mapsto(y,x)$ of the plane. 
For every $i\in[k]$ it induces a~bijection between 
$$
C_i(n)\times C_{L(i,\,k)}(n)\;\text{ and }C_{3k-i+1}(n)\times C_{L(3k-i+1,\,k)}(n)\;,
$$
and for every $i\in[k+1,2k]_{\Z}$ a~bijection between 
$$
C_i(n)\times C_{I_2(i,\,k)}(n) \;\text{ and }C_{3k-i+1}(n)\times C_{K_2(3k-i+1,\,k)}(n)\;.
$$

\noindent
{\bf 3.~Expressing $S_1(n)$, $S_2'(n)$ and $S_2''(n)$ by sine and cosine in 
(\ref{S1asympt})--(\ref{S2''explic}). }We substitute for each cardinality 
$|\cdots|$ on the right-hand sides of equations (\ref{cardinal}) with $j=1$, (\ref{S2'}) 
and (\ref{S2''}) the estimate 
of Proposition~\ref{prop_aproximacni1}. To this end we introduce some notation. 
Recall that for $i\in[3k]$,  $i\ga=i\pi/6k$. We consider the functions
$$
f_1(t)=\sin t-\cos t,\ f_2(t)=-\sin t-\cos t,\ f_3(t)=-\sin t+\cos t 
$$
and 
$$
f_4(t)=\sin t+\cos\;;
$$
we use $f_j(t)$ in the $j$-th quadrant. We label quadrants counter-clockwisely, 
starting with the first one where $x,y>0$. 
For any function $f\colon[0,2\pi)=C\to\R$ and any interval 
$I=[a,b]_{\Z}\sus[12k]$, $a,b\in[12k]$, we write 
$$
[f]_I=\big[f\big]_{(a-1)\ga}^{b\ga}\;.
$$
Note that if $I=[a,b]_{\Z}$ is any of the above intervals $I_j,J_j$ and $K_2$, $j\in[3]$, then always $a-1\le b$.

Thus using Proposition~\ref{prop_aproximacni1}, for any fixed $\ep>0$ and all $k,n\in\N$ 
satisfying the basic assumption we get that
\begin{eqnarray}
    S_1(n)&=&S_1\cdot n^2+O_{\ep}(kn^{1+\ep})\;,\label{S1asympt}\\
    S_2'(n)&=&S_2'\cdot n^2+O_{\ep}(kn^{1+\ep})\;\text{ and }\label{S2'asympt}\\
    S_2''(n)&=&S_2''\cdot n^2+O_{\ep}(kn^{1+\ep})\;,\label{S2''asympt}
    \;\text{ where}\\
    S_1&=&\sum_{i=1}^{k} [f_1(t)]_{(i-1)\ga}^{i\ga}\cdot\big([f_2(t)]_{I_1}+[f_3(t)]_{J_1}\big)\;,\label{S1explic}\\
    S_2'&=&\sum_{i=k+1}^{2k} [f_1(t)]_{(i-1)\ga}^{i\ga}\cdot[f_2(t)]_{I_2}\label{S2'explic}\;\text{ and}\\
    S_2''&=&\sum_{i=k+1}^{2k} [f_1(t)]_{(i-1)\ga}^{i\ga}\cdot[f_3(t)]_{J_2}\;.\label{S2''explic}
\end{eqnarray}

\noindent
{\bf 4.~Computing the sum $S_1$ in (\ref{S1_vysledek}). }We evaluate sums $S_1$, $S_2'$ 
and $S_2''$ in equations (\ref{S1explic}), (\ref{S2'explic}) and (\ref{S2''explic}), 
respectively. We begin with sum (\ref{S1explic}). For any real $x\ge 0$ we define ($i\in\N$)
\begin{equation}\label{Aofx}
A(x)=\sum_{i\le x}[f_1(t)]_{(i-1)\ga}^{i\ga}=f_1(\lfloor x\rfloor\ga)+1\;,
\end{equation}
with the empty sum defined as $0$. One has that
\begin{equation}\label{f2f3}
[f_2(t)]_{I_1(i,\,k)}+[f_3(t)]_{J_1(i,\,k)}=2+g(t)_{t=i}
\end{equation}
where
$$
{\textstyle
g(t)=\sin(\frac{2\pi}{3}+\ga t)+\cos(\frac{2\pi}{3}+\ga t)-
\sin(\frac{4\pi}{3}+\ga t-\ga)+\cos(\frac{4\pi}{3}+\ga t-\ga)\;.
}
$$
Using in equation (\ref{S1explic}) expressions (\ref{Aofx}) and (\ref{f2f3}), we get by Proposition~\ref{prop_sumAbel}, invoked on the second line, that 
\begin{eqnarray}
S_1&=&\sum_{i=1}^k[f_1(t)]_{(i-1)\ga}^{i\ga}\cdot(2+g(i))=2[f_1(t)]_0^{\pi/6}+\sum_{i=1}^k [f_1(t)]_{(i-1)\ga}^{i\ga}\cdot g(i)\notag\\
&=&3-\sqrt{3}+[A(t)g(t)]_0^k-\int_0^k A(t)\cdot g'(t)\;\mathrm{dt}\notag\\
&=&3-\sqrt{3}+[f_1(\lfloor t\rfloor\ga)\cdot g(t)]_0^k-\int_0^k f_1(\lfloor t\rfloor\ga)\cdot g'(t)\;\mathrm{dt}\notag\\
&=&
3-\sqrt{3}+\frac{1}{2}-\int_0^k f_1(\lfloor t\rfloor\ga)\cdot g'(t)\;\mathrm{dt}\notag\\
&=&
3-\sqrt{3}+\frac{1}{2}-\int_{2\pi/3}^{5\pi/6}f_1(\lfloor (u-2\pi/3)/\ga\rfloor\ga)\cdot(\cos u-\sin u)\;\mathrm{du}+\notag\\
&&+\;\int_{4\pi/3}^{3\pi/2}f_1(\lfloor (u-4\pi/3)/\ga\rfloor\ga)\cdot(\cos u+\sin u)\;\mathrm{du}\;.\label{dva_integr}
\end{eqnarray}
We evaluate the last two integrals with precision $O(1/k)$. We use the formulas $\sin(x+y)=\sin x\cos y+\cos x\sin y$, $\cos(x+y)=\cos x\cos y-\sin x\sin y$, $\sin(x+O(1/k))=\sin x+O(1/k)$, $\cos(x+O(1/k))=\cos x+O(1/k)$ (recall that $\gamma=\pi/6k$) and get that
\begin{eqnarray}
&&f_1(\lfloor(u-2\pi/3)/\ga\rfloor\ga)=f_1(u-2\pi/3-O(1/k))\notag\\
&&{\textstyle
=\sin u\cdot(-\frac{1}{2})+\cos u\cdot(-\frac{\sqrt{3}}{2})-\cos u\cdot(-\frac{1}{2})+\sin u\cdot(-\frac{\sqrt{3}}{2})+O(1/k)\notag}\\
&&={\textstyle
-\frac{1+\sqrt{3}}{2}\sin u+\frac{1-\sqrt{3}}{2}\cos u+O(1/k)
}\label{f1Delta1}
\end{eqnarray}
and
\begin{eqnarray}
&&f_1(\lfloor(u-4\pi/3)/\ga\rfloor\ga)=f_1(u-4\pi/3-O(1/k))\notag\\
&&{\textstyle
=\sin u\cdot(-\frac{1}{2})+\cos u\cdot\frac{\sqrt{3}}{2}-
\cos u\cdot(-\frac{1}{2})+\sin u\cdot\frac{\sqrt{3}}{2}+O(1/k)\notag}\\
&&{\textstyle=\frac{\sqrt{3}-1}{2}\sin u+\frac{1+\sqrt{3}}{2}\cos u+O(1/k)\;.}\label{f1Delta2}
\end{eqnarray}
After substituting expansions (\ref{f1Delta1}) and  (\ref{f1Delta2}) 
in the respective integrands in equation (\ref{dva_integr}) we see that the first integral in equation (\ref{dva_integr}) is
$$
{\textstyle-\int\sin u\cos u\;\mathrm{du}+\frac{1+\sqrt{3}}{2}\int\sin^2 u\;\mathrm{du}+\frac{1-\sqrt{3}}{2}\int\cos^2 u\;\mathrm{du}+O(1/k)}\;,
$$
with $\int=\int_{2\pi/3}^{5\pi/6}$,
and the second one is
$$
{\textstyle\sqrt{3}\int\sin u\cos u\;\mathrm{du}+\frac{\sqrt{3}-1}{2}\int\sin^2 u\;\mathrm{du}+\frac{1+\sqrt{3}}{2}\int\cos^2 u\;\mathrm{du}+O(1/k)}\;,
$$
with $\int=\int_{4\pi/3}^{3\pi/2}$. The former expression evaluates to (recall that $\sin(2x)=2\sin x\cdot\cos x$ and $\cos(2x)=\cos^2 x -\sin^2 x=1-2\sin^2 x=2\cos^2 x-1$)
\begin{eqnarray}
&&{\textstyle
\frac{1}{4}[\cos(2u)]_{2\pi/3}^{5\pi/6}+\frac{1+\sqrt{3}}{8}[2u-\sin(2u)]_{2\pi/3}^{5\pi/6}+\frac{1-\sqrt{3}}{8}[2u+\sin(2u)]_{2\pi/3}^{5\pi/6}+}\notag\\
&&+\;O(1/k)\notag\\
&&{\textstyle=\frac{1}{4}\cdot1
+\frac{1+\sqrt{3}}{8}\cdot\frac{\pi}{3}+\frac{1-\sqrt{3}}{8}\cdot\frac{\pi}{3}+O(1/k)\notag}\\
&&\label{integ1}{\textstyle=\frac{1+\pi/3}{4}+O(1/k)}\;,
\end{eqnarray}
and the latter to 
\begin{eqnarray}
&&{\textstyle
\frac{\sqrt{3}}{4}[-\cos(2u)]_{4\pi/3}^{3\pi/2}+\frac{\sqrt{3}-1}{8}[2u-\sin(2u)]_{4\pi/3}^{3\pi/2}+\frac{1+\sqrt{3}}{8}[2u+\sin(2u)]_{4\pi/3}^{3\pi/2}+}\notag\\
&&+\;O(1/k)\notag\\
&&{\textstyle=\frac{\sqrt{3}}{8}
+\frac{\sqrt{3}-1}{8}\cdot\left(\frac{\pi}{3}+\frac{\sqrt{3}}{2}\right)+\frac{1+\sqrt{3}}{8}\cdot\left(\frac{\pi}{3}-\frac{\sqrt{3}}{2}\right)+O(1/k)\notag}\\
&&{\textstyle=\frac{\pi\sqrt{3}}{12}+O(1/k)}\;.\label{integ2}
\end{eqnarray}
Substituting values (\ref{integ1}) and (\ref{integ2}) for the first and second integral in equation  (\ref{dva_integr}), respectively, 
we therefore get that
\begin{eqnarray}
    S_1&=&3-\sqrt{3}+\frac{1}{2}-\frac{1+\pi/3}{4}+\frac{\pi\sqrt{3}}{12}+O(1/k)\notag\\
    &=&\frac{39-12\sqrt{3}+\pi(\sqrt{3}-1)}{12}
    +O(1/k)\;.\label{S1_vysledek}
\end{eqnarray}

\noindent
{\bf 5.~Computing the sum $S_2'$ in (\ref{S2'_vysledek}). }We evaluate the second sum (\ref{S2'explic}). Now
\begin{equation}\label{f2f4}
[f_2(t)]_{I_2(i,\,k)}=:1+h(t)_{t=i}
\end{equation}
where
$$
h(t)=\sin(2\pi/3+\ga t)+\cos(2\pi/3+\ga t)\;.
$$
Using in equation (\ref{S2'explic}) expressions (\ref{Aofx}) and (\ref{f2f4}), we get by Proposition~\ref{prop_sumAbel}, invoked on the second line, that 
\begin{eqnarray}
S_2'&=&\sum_{i=k+1}^{2k}[f_1(t)]_{(i-1)\ga}^{i\ga}\cdot(1+h(i))=[f_1(t)]_{\pi/6}^{\pi/3}+\sum_{i=k+1}^{2k} [f_1(t)]_{(i-1)\ga}^{i\ga}\cdot h(i)\notag\\
&=&\sqrt{3}-1+[A(t)h(t)]_k^{2k}-\int_k^{2k} A(t)\cdot h'(t)\;\mathrm{dt}\notag\\
&=&\sqrt{3}-1+[f_1(\lfloor t\rfloor\ga)\cdot h(t)]_k^{2k}-\int_k^{2k} f_1(\lfloor t\rfloor\ga)\cdot h'(t)\;\mathrm{dt}\notag\\
&=&
\sqrt{3}-\frac{3}{2}-\int_k^{2k} f_1(\lfloor t\rfloor\ga)\cdot h'(t)\;\mathrm{dt}\notag\\
&=&
\sqrt{3}-\frac{3}{2}-\int_{5\pi/6}^{\pi}f_1(\lfloor (u-2\pi/3)/\ga\rfloor\ga)\cdot(\cos u-\sin u)\;\mathrm{du}\;.\label{jeden_integr}
\end{eqnarray}
We substitute expansion (\ref{f1Delta1})
in the integrand in equation (\ref{jeden_integr}) and see that the integral in equation (\ref{jeden_integr}) is again
$$
{\textstyle-\int\sin u\cos u\;\mathrm{du}+\frac{1+\sqrt{3}}{2}\int\sin^2 u\;\mathrm{du}+\frac{1-\sqrt{3}}{2}\int\cos^2 u\;\mathrm{du}+O(1/k)}
$$
but now $\int=\int_{5\pi/6}^{\pi}$. Now this expression evaluates to
\begin{eqnarray}
&&{\textstyle
\frac{1}{4}[\cos(2u)]_{5\pi/6}^{\pi}+\frac{1+\sqrt{3}}{8}[2u-\sin(2u)]_{5\pi/6}^{\pi}+\frac{1-\sqrt{3}}{8}[2u+\sin(2u)]_{5\pi/6}^{\pi}+}\notag\\
&&+\;O(1/k)\notag\\
&&{\textstyle=\frac{1}{4}\cdot\frac{1}{2}
+\frac{1+\sqrt{3}}{8}\cdot\left(\frac{\pi}{3}-\frac{\sqrt{3}}{2}\right)+\frac{1-\sqrt{3}}{8}\cdot\left(\frac{\pi}{3}+\frac{\sqrt{3}}{2}\right)+O(1/k)\notag}\\
&&\label{jeden_integr2}{\textstyle=\frac{\pi}{12}-\frac{1}{4}+O(1/k)}\;.
\end{eqnarray}
Substituting value (\ref{jeden_integr2}) for the integral in equation  (\ref{jeden_integr}), 
we get that
\begin{equation}\label{S2'_vysledek}
    {\textstyle
S_2'=\sqrt{3}-\frac{3}{2}-\left(\frac{\pi}{12}-\frac{1}{4}\right)+O(1/k)=\sqrt{3}-\frac{5}{4}-\frac{\pi}{12}+O(1/k)\;.}
\end{equation}

\noindent
{\bf 6.~Computing the sum $S_2''$ in (\ref{S2''_vysledek}). }We proceed to the third and last sum (\ref{S2''explic}).
Since the interval $J_2(i,k)$ is independent of $i$, the sum $S_2''$ can be computed easily and precisely:
\begin{eqnarray}\label{S2scarou}
S_2''&=&\sum_{i=k+1}^{2k} [f_1(t)]_{(i-1)\ga}^{i\ga}\cdot[f_3(t)]_{J_2(i,\,k)}\notag\\
&=&[\sin t-\cos t]_{\pi/6}^{\pi/3}\cdot[-\sin t+\cos t]_{\pi}^{3\pi/2}\notag\\
&=&2(\sqrt{3}-1)\;.\label{S2''_vysledek}
\end{eqnarray}

\noindent
{\bf 7. Conclusion of the computation. }Using equations (\ref{DnaD1qn}), (\ref{Ch1qaSj}), (\ref{S2_sumofthree}), (\ref{symmetry}), (\ref{S1asympt}), (\ref{S2'asympt}), (\ref{S2''asympt}), (\ref{S1_vysledek}), (\ref{S2'_vysledek}), (\ref{S2''_vysledek}) and setting $k=k(n)=\lfloor\sqrt{n}\rfloor$, we deduce that for any fixed $\ep>0$ and all large $n\in\N$, 
\begin{eqnarray}
  |\mathrm{Ch}(n)|&=&{\textstyle 4\cdot\big(2\cdot\frac{39-12\sqrt{3}+\pi(\sqrt{3}-1)}{12}+2\cdot(\sqrt{3}-\frac{5}{4}-\frac{\pi}{12})+2(\sqrt{3}-1)\big)\cdot n^2+\notag}\\
  &&+\;O_{\ep}(n^{3/2+\ep})\notag\\
  &=&\frac{24+24\sqrt{3}+\pi(2\sqrt{3}-4)}{3}\cdot n^2+O_{\ep}(n^{3/2+\ep})\;.
  \label{asympt_Ch}
\end{eqnarray}
Finally, by equation (\ref{asympt_Ch}) and Corollary~\ref{cor_sol3} we see that
\begin{eqnarray*}
  \mathrm{Pr_{Be,\,\Box}}=\lim_{n\to\infty}\frac{|\mathrm{Ch}(n)|}{|C(n)|^2}&=&\frac{24(1+\sqrt{3})+\pi(2\sqrt{3}-4)}{3\cdot 64}\\
  &=&\frac{1+\sqrt{3}}{8}-\frac{\pi(2-\sqrt{3})}{96}=0.33273\ds\;.
\end{eqnarray*}
\eproof

\noindent
This is close to but different from Bertrand's Solution~1 $\mathrm{Pr_{Be}}=\frac{1}{3}$.

\section{Suggestions of a~referee and concluding remarks}

A~referee, dismissive of the approach here, wrote the following.

\bigskip\noindent
 $(\ds)$\\
``It is is easy to see that in the limit $n\to\infty$ this gives rise to the following distribution 
on the circle: with the usual parametrization $\varphi\mapsto e^{i\varphi}$, $\varphi\in[0,2\pi)$ 
of the unit circle, the probability density function is
$$
f(\varphi)=\frac{|\sin(\varphi)|+|\cos(\varphi)|}{8}.
$$
In the limit problem, one needs to choose the two endpoints of the chord independently, 
each according to this distribution. So the probability in question is
$$
\int_0^{2\pi}\int_{2\pi/3}^{4\pi/3}f(x)f(x+y)\,\mathrm{d}x\mathrm{d}y.
$$
It is a~straightforward (although not pleasent [sic])  task to compute this integral.''\\
$(\ds)$\\
``I did not check the author's computations (that take up roughly 12 pages) but I am certain 
that they are unnecessarily complicated and the result could be obtained in a~few pages 
(by proving a~weak convergence for the discrete approximations and computing the above integral 
for the limit problem).''\\
 $(\ds)$

\bigskip\noindent
In the sentence ``It is is easy to see that $\ds$'' we should get a~credit and Proposition~\ref{prop_aproximacni1} 
should be mentioned because the density $f(\varphi)$ follows 
from the asymptotics in Proposition~\ref{prop_aproximacni1}, not the other way around (as 
is much more common in probabilistic computations). Were the values 
of the arithmetic function $r_2(m)$\,---\,see the proof of Proposition~\ref{prop_onExceptionals}\,---\,not small 
enough, $f(\varphi)$ would be different, if it existed at all. The obvious bound $r_2(m)=O(m^{1/2})$
does not suffice, nor even $r_2(m)=O(m^{1/4})$. In reality it is not 
true that ``It is is easy to see that $\ds$''. 

MAPLE (2022 edition) cannot compute the above double integral in a~closed form, instead it expresses it, 
by computing the 
inner integral, in the form
$$
I:=\int_0^{2\pi}I(x)\,\mathrm{d}x
$$
with a~complicated integrand $I(x)$ defined by considering various cases. More precisely, this result was only 
obtained for some forms of the integrand $f(x)f(x+y)$. For other trivial variations, MAPLE could not compute 
anything or was running inconclusively for one hour before it was terminated. It remains to be seen if this reflects some 
inherent complexity of the integral or if it results from imperfection of algorithms in MAPLE. 

Fortunately, MAPLE can 
compute approximations of the integral, 
$$
I=0.33273773412864161039\ds\;.
$$
These decimal digits agree with our value of $\mathrm{Pr_{Be,\Box}}$. We have therefore a~non-rigorous 
``proof'' of the fact that
$$
I=\frac{1+\sqrt{3}}{8}-\frac{\pi(2-\sqrt{3})}{96}\;.
$$
To make it rigorous would mean (i) to prove ``a~weak convergence for the discrete approximations'' or (ii) to 
compute the double integral in closed form. We are less optimistic about the simplicity of a~proof in (i) than the referee but this can be decided only by providing a~detailed and rigorous proof in (i). 
We thank the referee for the suggested alternative approach to asymptotics of sums counting long 
chords or similar quantities.

Returning to our view of Bertrand's paradox, we suggest to restate the problem and, especially, what it means 
to solve it in the following terms. One should not be searching, like some analyst of old trying to find the true sum of the series $1-1+1-1+1-1+\ds$, for the true probability of long chords in the unit geometric circle. 
\begin{quote}
Instead, to resolve Bertrand's paradox, one should devise a~framework that preserves
as much as possible from the spirit of the original intuitive problem (suppose we throw erratically, ``randomly'' 
sticks on a~circle drawn by chalk on the sidewalk, what is the chance that $\ds$) and at the same time makes the solution unique.
\end{quote}
The approximation by boxes in this article is a~possible solution of this metaproblem. Of course,
one cannot get rid of non-uniqueness completely, it is only shifted to the higher level.

To close, as for some ideas for further research we only mention treating the Buffon needle 
problem in our framework  (\cite{klaz_buff}). 
The method of discretization is not new in geometric probability\,---\,this author learned it from his colleague P.~Valtr in \cite{valt1,valt2,valt3}, but its application here to Bertrand's paradox is new.

\medskip\noindent
{\em Department of Applied Mathematics\\
Faculty of Mathematics and Physics\\
Charles University\\
Malostransk\'e n\'am. 25\\
118 00 Praha\\
Czechia}

\end{document}